\documentclass[a4paper,10pt,12pt]{article}
\usepackage{amsmath, amssymb, latexsym, amscd, amsthm,amsfonts,amstext}
\usepackage[mathscr]{eucal}
\usepackage{graphicx}
\usepackage{subfig}

\numberwithin{equation}{section}
\input{tcilatex}
\makeatletter
\renewcommand\tableofcontents{    \@starttoc{toc}}
\makeatother

\begin{document}

\title{A universal regularization method for ill-posed Cauchy problems for
quasilinear partial differential equations}
\author{Michael V. Klibanov \\
\\
Department of Mathematics \& Statistics, University of North Carolina \\
at Charlotte, Charlotte, NC 28223, USA\\
Email: \texttt{\ }mklibanv@uncc.edu }
\date{}
\maketitle

\begin{abstract}
For the first time, a globally convergent numerical method is presented for
ill-posed Cauchy problems for quasilinear PDEs. The key idea is to use
Carleman Weight Functions to construct globally strictly convex
Tikhonov-like cost functionals.
\end{abstract}



\section{Introduction}

\label{sec:1}

This is the first publication where a globally convergent numerical method
is presented for ill-posed Cauchy problems for a broad class of quasilinear
Partial Differential Equations (PDEs). All previous works were concerned
only with the linear case. First, we present the general framework of this
method. Next, we specify that framework for ill-posed Cauchy problems for
quasilinear elliptic, parabolic and hyperbolic PDEs.

Let $G\subset \mathbb{R}^{n}$ be a bounded domain with a piecewise smooth
boundary $\partial G.$ Let $\Gamma \subseteq \partial G$ be a part of the
boundary.\ Let $A$ be a quasilinear Partial Differential Operator (PDO) of
the second order acting on functions $u$ defined in the domain $G$ (details
are given in section 2). Suppose that $\Gamma $ is not a characteristic
hypersurface of $A$. Consider the Cauchy problem for the quasilinear PDE
generated by the operator $A$ with the Cauchy data at $\Gamma .$ Suppose
that this problem is ill-posed in the classical sense, e.g. the Cauchy
problem for a quasilinear elliptic equation. In this paper we construct a
universal globally convergent regularization method for solving such
problems. This method works for those PDOs, for which Carleman estimates are
valid for linearized versions of those operators. Since these estimates are
valid for three main types of linear PDOs of the second order, elliptic,
parabolic and hyperbolic ones, then our method is a quite general one.

The key idea of this paper is to construct a weighted Tikhonov-like
functional.\ The weight is the Carleman Weight Function (CWF), i.e. the
function which is involved in the Carleman estimate for the principal part
of the corresponding PDO. Given a ball $B\left( R\right) $ of an arbitrary
radius $R$ in a certain Sobolev space, one can choose the parameter $\lambda
\left( R\right) $ of this CWF so large that the weighted Tikhonov-like
functional is strictly convex on $B\left( R\right) $ for all $\lambda \geq
\lambda \left( R\right) $. The strict convexity, in turn guarantees the
convergence of the gradient method of the minimization of this functional if
starting from an arbitrary point of the ball $B\left( R\right) $. Since
restrictions of the radius $R$ are not imposed, then this is \emph{global
convergence}. On the other hand, the major problem with conventional
Tikhonov functionals for nonlinear ill-posed problems is that they usually
have many local minima and ravines. The latter implies convergence of the
gradient method only if it starts in a sufficiently small neighborhood of
the solution, i.e. the local convergence.

First globally strictly convex cost functionals for nonlinear ill-posed
problems were constructed by Klibanov \cite{Klib97,Kpar} for Coefficient
Inverse Problems (CIPs), using CWFs. Based on a modification of this idea,
some numerical studies for 1-d CIPs were performed in the book of Klibanov
and Timonov \cite{KT}. Recently there is a renewed interest to this topic,
see Baudoin, DeBuhan and Ervedoza \cite{Baud}, Beilina and Klibanov \cite%
{BKconv}, Klibanov and Thanh \cite{KNT} and Klibanov and Kamburg \cite{KK}.
In particular, the paper \cite{KNT} contains some numerical results.
However, globally strictly convex cost functionals were not constructed for
ill-posed Cauchy problems for quasilinear PDEs in the past.

Concerning applications of Carleman estimates to inverse problems, we refer
to the method, which was first proposed by Bukhgeim and Klibanov \cite%
{BukhK,Bukh1,Klib1}. The method of \cite{BukhK,Bukh1,Klib1} was originally
designed for proofs of uniqueness theorems for CIPs with single measurement
data, see, e.g. some follow up works of Bukhgeim \cite{Bukh}, Klibanov \cite%
{K92}, Klibanov and Timonov \cite{KT}, surveys of Klibanov \cite{Ksurvey}
and Yamamoto \cite{Y}, as well as sections 1.10 and 1.11 of the book of
Beilina and Klibanov \cite{BK}.

There is a huge number of publications about ill-posed Cauchy problems for
linear PDEs. Hence, we now refer only to a few of them. We note first that
in the conventional case of a linear ill-posed problem, Tikhonov
regularization functional is generated by a bounded linear operator, see,
e.g. books of Ivanov, Vasin and Tanana \cite{Iv} and Lavrentiev, Romanov and
Shishatskii \cite{LRS}. Unlike this, the idea, which is the closest one to
this paper, is to use unbounded linear PDOs as generators of Tikhonov
functionals for ill-posed Cauchy problems. This idea was first proposed in
the book of Lattes and Lions \cite{LL}. In \cite{LL} Tikhonov functionals
were written in the variational form. Since Carleman estimates were not used
in \cite{LL}, then convergence rates for minimizers were not established.
The idea of using Carleman estimates for establishing convergence rates of
minimizers of those Tikhonov functionals was first proposed in works of
Klibanov and Santosa \cite{KS} and of Klibanov and Malinsky \cite{KM}. Next,
it was explored in works of Klibanov with coauthors \cite{Cao,ClK,KKKN},
where\ accurate numerical results were obtained; also see a recent survey 
\cite{Klib}. In addition, Bourgeois \cite{B5,B9} and Bourgeois and Dard\'{e}
have used this idea for the case of the Cauchy problem for the Laplace
equation, see, e.g. \cite{B5,B9}. We also refer to Berntsson, Kozlov,
Mpinganzima and Turesson \cite{Kozlov3}, Hao and Lesnic \cite{HaoLes},
Kabanikhin \cite{Kab}, Kabanikhin and Karchevsky \cite{Karch1}, Kozlov,
Maz'ya and Fomin \cite{Kozlov1} and Li, Xie and Zou \cite{Zou} for some
other numerical methods for ill-posed Cauchy problems for linear PDEs.

All functions considered below are real valued ones. In section 2 we present
the general framework of our method. Next, in sections 3, 4 and 5 we specify
this framework for quasilinear elliptic, parabolic and hyperbolic PDEs
respectively.

\section{The method}

\label{sec:2}

Let $A$ be a quasilinear PDO of the second order in $G$ with its principal
part $A_{0},$ 
\begin{eqnarray}
A\left( u\right) &=&\sum\limits_{\left\vert \alpha \right\vert =2}a_{\alpha
}\left( x\right) D^{\alpha }u+A_{1}\left( x,\nabla u,u\right) ,
\label{2.100} \\
A_{0}u &=&\sum\limits_{\left\vert \alpha \right\vert =2}a_{\alpha }\left(
x\right) D^{\alpha }u,  \label{2.1001} \\
\text{ where functions }a_{\alpha } &\in &C^{1}\left( \overline{G}\right)
,A_{1}\in C\left( \overline{G}\right) \times C^{3}\left( \mathbb{R}%
^{n+1}\right) .  \label{2.1002}
\end{eqnarray}%
Let $k_{n}=\left[ n/2\right] +2,$ where $\left[ n/2\right] $ is the largest
integer which does not exceed the number $n/2.$ By the embedding theorem 
\begin{equation}
H^{k_{n}}\left( G\right) \subset C^{1}\left( \overline{G}\right) \text{ and }%
\left\Vert f\right\Vert _{C^{1}\left( \overline{G}\right) }\leq C\left\Vert
f\right\Vert _{H^{k_{n}}\left( G\right) },\forall f\in H^{k_{n}}\left(
G\right) ,  \label{2.10005}
\end{equation}%
where the constant $C=C\left( n,G\right) >0$ depends only on listed
parameters.

\textbf{Cauchy Problem}. \emph{Let the hypersurface }$\Gamma \subseteq
\partial G$\emph{\ and assume that }$\Gamma $ \emph{is not a characteristic
hypersurface of the operator }$A_{0}.$\emph{\ Consider the following Cauchy
problem for the operator }$A$\emph{\ with the Cauchy data }$g_{0}\left(
x\right) ,g_{1}\left( x\right) ,$\emph{\ }%
\begin{equation}
A\left( u\right) =0,  \label{2.1003}
\end{equation}%
\begin{equation}
u\mid _{\Gamma }=g_{0}\left( x\right) ,\partial _{n}u\mid _{\Gamma
}=g_{1}\left( x\right) ,x\in \Gamma .  \label{2.1004}
\end{equation}%
\emph{Determine the solution }$u\in H^{k_{n}}\left( G\right) $\emph{\ of the
problem (\ref{2.1003}), (\ref{2.1004}) either in the entire domain }$G$\emph{%
\ or at least in its subdomain.}

\subsection{The Carleman estimate}

\label{sec:2.1}

Let the function $\xi \in C^{2}\left( \overline{G}\right) $ and $\left\vert
\nabla \xi \right\vert \neq 0$ in $\overline{G}.$ For a number $c\geq 0$
denote 
\begin{equation}
\xi _{c}=\left\{ x\in \overline{G}:\xi \left( x\right) =c\right\}
,G_{c}=\left\{ x\in G:\xi \left( x\right) >c\right\} .  \label{2.0}
\end{equation}%
We assume that $G_{c}\neq \varnothing .$ Choose a sufficiently small $%
\varepsilon >0$ such that $G_{c+2\varepsilon }\neq \varnothing .$ Obviously, 
$G_{c+2\varepsilon }\subset G_{c+\varepsilon }\subset G_{c}.$ Let $\Gamma
_{c}=\left\{ x\in \Gamma :\xi \left( x\right) >c\right\} \neq \varnothing .$
Hence, the boundary of the domain $G_{c}$ consists of two parts,%
\begin{equation}
\partial G_{c}=\partial _{1}G_{c}\cup \partial _{2}G_{c},\partial
_{1}G_{c}=\xi _{c},\partial _{2}G_{c}=\Gamma _{c}.  \label{2.1}
\end{equation}

Let $\lambda >1$ be a large parameter. Consider the function $\varphi
_{\lambda }\left( x\right) ,$%
\begin{equation}
\varphi _{\lambda }\left( x\right) =\exp \left( \lambda \xi \left( x\right)
\right) .  \label{2.2}
\end{equation}%
It follows from (\ref{2.1}) and (\ref{2.2}) that 
\begin{equation}
\min_{\overline{G}_{c}}\varphi _{\lambda }\left( x\right) =\varphi _{\lambda
}\left( x\right) \mid _{\xi _{c}}=e^{\lambda c}.  \label{2.3}
\end{equation}

\textbf{Definition 2.1}. \emph{We say that the operator }$A_{0}$\emph{\
admits the pointwise\ Carleman estimate in the domain }$G_{c}$\emph{\ with
the CWF }$\varphi _{\lambda }\left( x\right) $\emph{\ if there exist
constants }$\lambda _{0}\left( G_{c},A_{0}\right) >1,C_{1}\left(
G_{c},A_{0}\right) >0$\emph{\ depending only on the domain }$G$\emph{\ and
the operator }$A_{0}$\emph{\ such that the following estimate holds}%
\begin{eqnarray}
\left( A_{0}u\right) ^{2}\varphi _{\lambda }^{2}\left( x\right) &\geq
&C_{1}\lambda \left( \nabla u\right) ^{2}\varphi _{\lambda }^{2}\left(
x\right) +C_{1}\lambda ^{3}u^{2}\varphi _{\lambda }^{2}\left( x\right) +%
\func{div}U,  \label{2.6} \\
\forall \lambda &\geq &\lambda _{0},\forall u\in C^{2}\left( \overline{G}%
\right) ,\forall x\in G_{c}.  \label{2.7}
\end{eqnarray}%
\emph{In (\ref{2.6}) the vector function }$U\left( x\right) $ \emph{%
satisfies the following estimate} 
\begin{equation}
\left\vert U\left( x\right) \right\vert \leq C_{1}\lambda ^{3}\left[ \left(
\nabla u\right) ^{2}+u^{2}\right] \varphi _{\lambda }^{2}\left( x\right)
,\forall x\in G_{c}.  \label{2.8}
\end{equation}

\subsection{The main result and the gradient method}

\label{sec:2.2}

Let $R>0$ be an arbitrary number. Denote 
\begin{equation}
B\left( R\right) =\left\{ u\in H^{k_{n}}\left( G_{c}\right) :\left\Vert
u\right\Vert _{H^{k_{n}}\left( G_{c}\right) }<R,u\mid _{\Gamma
_{c}}=g_{0}\left( x\right) ,\partial _{n}u\mid _{\Gamma _{c}}=g_{1}\left(
x\right) \right\} ,  \label{1}
\end{equation}%
\begin{equation}
H_{0}^{k_{n}}\left( G_{c}\right) =\left\{ u\in H^{k_{n}}\left( G_{c}\right)
:u\mid _{\Gamma _{c}}=\partial _{n}u\mid _{\Gamma _{c}}=0\right\} .
\label{2}
\end{equation}%
To solve the above Cauchy problem, we consider the following minimization
problem.

\textbf{Minimization Problem}. \emph{Assume that the operator }$A_{0}$\emph{%
\ satisfies the Carleman estimate (\ref{2.6}), (\ref{2.7}) for a certain
number }$c\geq 0.$ \emph{Minimize the functional }$J_{\lambda ,\alpha
}\left( u\right) $ \emph{in} \emph{(\ref{2.1005})} \emph{on the set }$%
B\left( R\right) $\emph{, where}%
\begin{equation}
J_{\lambda ,\beta }\left( u\right) =e^{-2\lambda \left( c+\varepsilon
\right) }\dint\limits_{G_{c}}\left[ A\left( u\right) \right] ^{2}\varphi
_{\lambda }^{2}dx+\beta \left\Vert u\right\Vert _{H^{k_{n}}\left(
G_{c}\right) }^{2},  \label{2.1005}
\end{equation}%
\emph{where }$\beta \in \left( 0,1\right) $\emph{\ is the regularization
parameter.}

Thus, by solving this problem we approximate the function $u$ in the
subdomain $G_{c}\subset G.$ The multiplier $e^{-2\lambda \left(
c+\varepsilon \right) }$ is introduced to balance two terms in the right
hand side (\ref{2.1005}): to allow to have $\beta \in \left( 0,1\right) .$
Theorem 2.1 is the main result of this paper. Note that since $e^{-\lambda
\varepsilon }<<1$ for sufficiently large $\lambda ,$ then the requirement of
this theorem $\beta \in \left( e^{-\lambda \varepsilon },1\right) $ means
that the regularization parameter $\beta $ can change from being very small
and up to unity.

\textbf{Theorem 2.1}. \emph{Let }$R>0$\emph{\ be an arbitrary number. Let }$%
B\left( R\right) $\emph{\ and }$H_{0}^{k_{n}}\left( G_{c}\right) $\emph{\ be
the sets defined in (\ref{1}) and (\ref{2}) respectively. Then for every
point }$u\in B\left( R\right) $\emph{\ there exists the Fr\'{e}chet
derivative }$J_{\lambda ,\alpha }^{\prime }\left( u\right) \in
H_{0}^{k_{n}}\left( G_{c}\right) .$\emph{\ Assume that the operator }$A_{0}$%
\emph{\ admits the pointwise\ Carleman estimate in the domain }$G_{c}$\emph{%
\ of Definition 2.1 and let }$\lambda _{0}\left( G_{c},A_{0}\right) >1$\emph{%
\ be the constant of this definition. Then there exists a sufficiently large
number }$\lambda _{1}=\lambda _{1}\left( G_{c},A,R\right) >\lambda _{0}$%
\emph{\ such that for all }$\lambda \geq \lambda _{1}$\emph{\ and for every }%
$\beta \in \left( e^{-\lambda \varepsilon },1\right) $\emph{\ the functional 
}$J_{\lambda ,\beta }\left( u\right) $\emph{\ is strictly convex on the set }%
$B\left( R\right) ,$%
\begin{equation*}
J_{\lambda ,\beta }\left( u_{2}\right) -J_{\lambda ,\beta }\left(
u_{1}\right) -J_{\lambda ,\beta }^{\prime }\left( u_{1}\right) \left(
u_{2}-u_{1}\right)
\end{equation*}%
\begin{equation*}
\geq C_{2}e^{2\lambda \varepsilon }\left\Vert u_{2}-u_{1}\right\Vert
_{H^{1}\left( G_{c+2\varepsilon }\right) }^{2}+\frac{\beta }{2}\left\Vert
u_{2}-u_{1}\right\Vert _{H^{k_{n}}\left( G_{c}\right) }^{2},\forall
u_{1},u_{2}\in B\left( R\right) .
\end{equation*}%
\emph{Here the number }$C_{2}=C_{2}\left( A,R,c\right) $\emph{\ depends only
on listed parameters. }

We now show that this theorem implies the global convergence of the gradient
method of the minimization of the functional (\ref{2.1005}) on the set $%
B\left( R\right) $. Consider an arbitrary function $u_{1}\in B\left(
R\right) ,$ which is our starting point for iterations of this method. Let
the step size of the gradient method be $\gamma >0$. For brevity, we do not
indicate the dependence of functions $u_{n}$ on parameters $\lambda ,\beta
,\gamma $. Consider the sequence $\left\{ u_{n}\right\} _{n=1}^{\infty }$ of
the gradient method, 
\begin{equation}
u_{n+1}=u_{n}-\gamma J_{\lambda ,\alpha }^{\prime }\left( u_{n}\right)
,n=1,2,...  \label{2.200}
\end{equation}

\textbf{Theorem 2.2}. \emph{Let }$\lambda _{1}$\emph{\ be the parameter of
Theorem 2.1. Choose a number }$\lambda \geq \lambda _{0}.$ \emph{Let }$\beta
\in \left( e^{-\lambda \varepsilon },1\right) .$ \emph{Assume that the
functional }$J_{\lambda ,\beta }$\emph{\ achieves its minimal value on the
set }$B\left( R\right) $\emph{\ at a point }$u_{\min }\in B\left( R\right) $%
\emph{.}\ \emph{\ Then such a point }$u_{\min }$ \emph{is unique. Consider
the sequence (\ref{2.200}), where }$u_{1}\in B\left( R\right) $ \emph{is an
arbitrary point}. \emph{Assume that }$\left\{ u_{n}\right\} _{n=1}^{\infty
}\subset B\left( R\right) .$ \emph{Then there exists a sufficiently small
number }$\gamma =\gamma \left( \lambda ,\beta ,B\left( R\right) \right) \in
\left( 0,1\right) $ \emph{and a number }$q=q\left( \gamma \right) \in \left(
0,1\right) $\emph{\ such that the sequence (\ref{2.200}) converges to the
point }$u_{\min },$%
\begin{equation*}
\left\Vert u_{n+1}-u_{\min }\right\Vert _{H^{k_{n}}\left( G_{c}\right) }\leq
q^{n}\left\Vert u_{0}-u_{\min }\right\Vert _{H^{k_{n}}\left( G_{c}\right) }.
\end{equation*}

Since the starting point $u_{1}$ of the sequence (\ref{2.200}) is an
arbitrary point of the set $B\left( R\right) $ and a restriction on $R$ is
not imposed, then Theorem 2.2 claims the \emph{global convergence} of the
gradient method. This is the opposite to its local convergence for
non-convex functionals. Since it was shown in \cite{BKconv} how a direct
analog of Theorem 2.2 follows from an analog of Theorem 2.1, although for a
different cost functional, we do not prove Theorem 2.2 here for brevity. We
note that it is likely that similar global convergence results can be proven
for other versions of the gradient method, e.g. the conjugate gradient
method.\ However, we are not doing this here for brevity.

\subsection{Proof of Theorem 2.1}

In this proof $C_{2}=C_{2}\left( A,R,c\right) >0$ denotes different numbers
depending only on listed parameters. Let $u_{1},u_{2}\in B\left( R\right) $
be two arbitrary functions and let $h=u_{2}-u_{1}.$ Then (\ref{2.1006})
implies that 
\begin{equation}
h\in H_{0}^{k_{n}}\left( G_{c}\right) .  \label{2.11}
\end{equation}%
Let 
\begin{equation}
D=\left( A\left( u_{2}\right) \right) ^{2}-\left( A\left( u_{1}\right)
\right) ^{2}.  \label{2.12}
\end{equation}%
We now put the expression for $D$ in a form, which is convenient for us.
First, we recall that the Lagrange formula implies%
\begin{equation*}
f\left( y+z\right) =f\left( y\right) +f^{\prime }\left( y\right) z+\frac{%
z^{2}}{2}f^{\prime \prime }\left( \eta \right) ,\forall y,z\in \mathbb{R}%
,\forall f\in C^{2}\left( \mathbb{R}\right) ,
\end{equation*}%
where $\eta =\eta \left( y,z\right) $ is a number located between numbers $y$
and $z$. By (\ref{2.10005})%
\begin{equation}
\left\Vert h\right\Vert _{C^{1}\left( \overline{G}_{c}\right) }=\left\Vert
u_{2}-u_{1}\right\Vert _{C^{1}\left( \overline{G}_{c}\right) }\leq 2CR.
\label{2.13}
\end{equation}%
Hence, using this formula, (\ref{2.1001}) and (\ref{2.13}), we obtain for
the operator $A_{1}$ 
\begin{equation*}
A_{1}\left( x,\nabla \left( u_{1}+h\right) ,u_{1}+h\right) =A_{1}\left(
x,\nabla u_{1},u_{1}\right)
\end{equation*}%
\begin{equation*}
+\dsum\limits_{i=1}^{n}\partial _{u_{x_{i}}}A_{1}\left( x,\nabla
u_{1},u_{1}\right) h_{x_{i}}+\partial _{u}A_{1}\left( x,\nabla
u_{1},u_{1}\right) h+F\left( x,\nabla u_{1},u_{1},h\right) ,
\end{equation*}%
where the function $F$ satisfies the following estimate%
\begin{equation}
\left\vert F\left( x,\nabla u_{1},u_{1},h\right) \right\vert \leq
C_{2}\left( \left( \nabla h\right) ^{2}+h^{2}\right) ,\forall x\in
G_{c},\forall u_{1}\in B\left( R\right) .  \label{2.14}
\end{equation}%
Hence, 
\begin{equation*}
A\left( u_{1}+h\right) =A_{0}\left( u_{1}+h\right) +A_{1}\left( x,\nabla
\left( u_{1}+h\right) ,u_{1}+h\right) =
\end{equation*}%
\begin{equation*}
A\left( u_{1}\right) +\left[ A_{0}\left( h\right)
+\dsum\limits_{i=1}^{n}\partial _{u_{x_{i}}}A_{1}\left( x,\nabla
u_{1},u_{1}\right) h_{x_{i}}+\partial _{u}A_{1}\left( x,\nabla
u_{1},u_{1}\right) h\right] +F\left( x,\nabla u_{1},u_{1},h\right) .
\end{equation*}%
Hence, by (\ref{2.12})%
\begin{equation*}
D=2A\left( u_{1}\right) \left[ A_{0}\left( h\right)
+\dsum\limits_{i=1}^{n}\partial _{u_{x_{i}}}A_{1}\left( x,\nabla
u_{1},u_{1}\right) h_{x_{i}}+\partial _{u}A_{1}\left( x,\nabla
u_{1},u_{1}\right) h\right]
\end{equation*}%
\begin{equation}
+\left[ A_{0}\left( h\right) +\dsum\limits_{i=1}^{n}\partial
_{u_{x_{i}}}A_{1}\left( x,\nabla u_{1},u_{1}\right) h_{x_{i}}+\partial
_{u}A_{1}\left( x,\nabla u_{1},u_{1}\right) h\right] ^{2}+F^{2}  \label{2.15}
\end{equation}%
\begin{equation*}
+2\left[ A_{0}\left( h\right) +\dsum\limits_{i=1}^{n}\partial
_{u_{x_{i}}}A_{1}\left( x,\nabla u_{1},u_{1}\right) h_{x_{i}}+\partial
_{u}A_{1}\left( x,\nabla u_{1},u_{1}\right) h\right] F.
\end{equation*}

The expression in the first line of (\ref{2.15}) is linear with respect to $%
h $. We denote it as $Q\left( u_{1}\right) \left( h\right) .$ Consider the
linear functional 
\begin{equation}
\widetilde{J}_{u_{1}}\left( h\right) =\dint\limits_{G_{c}}Q\left(
u_{1}\right) \left( h\right) \varphi _{\lambda }^{2}dx+2\beta \left[ u_{1},h%
\right] ,  \label{2.16}
\end{equation}%
where $\left[ ,\right] $ denotes the scalar product in $H^{k_{n}}\left(
G_{c}\right) .$ Clearly, $\widetilde{J}_{u_{1}}\left( h\right)
:H_{0}^{k_{n}}\left( G_{c}\right) \rightarrow \mathbb{R}$ is a bounded
linear functional. Hence, by the Riesz theorem, there exists a single
element $P\left( u_{1}\right) \in H_{0}^{k_{n}}\left( G_{c}\right) $ such
that $\widetilde{J}_{u_{1}}\left( h\right) =\left[ P\left( u_{1}\right) ,h%
\right] ,\forall h\in H_{0}^{k_{n}}\left( G_{c}\right) .$ Furthermore, $%
\left\Vert P\left( u_{1}\right) \right\Vert _{H^{k_{n}}\left( G_{c}\right)
}=\left\Vert \widetilde{J}_{u_{1}}\right\Vert .$ This proves the existence
of the Frech\'{e}t derivative 
\begin{equation}
J_{\lambda ,\beta }^{^{\prime }}\left( u_{1}\right) =P\left( u_{1}\right)
\in H_{0}^{k_{n}}\left( G_{c}\right) .  \label{2.17}
\end{equation}

Let 
\begin{equation*}
S\left( x,u_{1},h\right) =D-2A\left( u_{1}\right) \left[ A_{0}\left(
h\right) +\dsum\limits_{i=1}^{n}\partial _{u_{x_{i}}}A_{1}\left( x,\nabla
u_{1},u_{1}\right) h_{x_{i}}+\partial _{u}A_{1}\left( x,\nabla
u_{1},u_{1}\right) h\right] .
\end{equation*}%
Then, using (\ref{2.13})-(\ref{2.15}) and the Cauchy-Schwarz inequality, we
obtain%
\begin{equation*}
S\geq \frac{1}{2}\left( A_{0}h\right) ^{2}-C_{2}\left( \left( \nabla
h\right) ^{2}+h^{2}\right) ,\forall x\in G_{c},\forall u_{1}\in B\left(
R\right) .
\end{equation*}%
Hence, using (\ref{2.16}) and (\ref{2.17}), we obtain%
\begin{equation*}
J_{\lambda ,\beta }\left( u_{1}+h\right) -J_{\lambda ,\beta }\left(
u_{1}\right) -J_{\lambda ,\beta }^{\prime }\left( u_{1}\right) \left(
h\right)
\end{equation*}%
\begin{equation}
\geq \frac{1}{2}e^{-2\lambda \left( c+\varepsilon \right)
}\dint\limits_{G_{c}}\left( A_{0}h\right) ^{2}\varphi _{\lambda
}^{2}dx-C_{2}e^{-2\lambda \left( c+\varepsilon \right)
}\dint\limits_{G_{c}}\left( \left( \nabla h\right) ^{2}+h^{2}\right) \varphi
_{\lambda }^{2}dx+\beta \left\Vert h\right\Vert _{H^{k_{n}}\left(
G_{c}\right) }^{2}.  \label{2.18}
\end{equation}

Next, integrate (\ref{2.6}) over the domain $G_{c},$ using the Gauss'
formula, (\ref{2.7}) and (\ref{2.8}). Next, replace $u$ with $h$ in the
resulting formula. Even though there is no guarantee that $h\in C^{2}\left( 
\overline{G}_{c}\right) ,$ still density arguments ensure that the resulting
inequality remains true. Hence, using (\ref{2.11}), we obtain 
\begin{equation*}
\frac{1}{2}e^{-2\lambda \left( c+\varepsilon \right)
}\dint\limits_{G_{c}}\left( A_{0}h\right) ^{2}\varphi _{\lambda }^{2}dx\geq 
\frac{C_{1}}{2}e^{-2\lambda \left( c+\varepsilon \right)
}\dint\limits_{G_{c}}\left( \lambda \left( \nabla h\right) ^{2}+\lambda
^{3}h^{2}\right) \varphi _{\lambda }^{2}dx
\end{equation*}%
\begin{equation*}
-\frac{C_{1}}{2}\lambda ^{3}e^{-2\lambda \varepsilon }\dint\limits_{\xi
_{c}}\left( \left( \nabla h\right) ^{2}+h^{2}\right) \varphi _{\lambda
}^{2}dx.
\end{equation*}%
Substituting this in (\ref{2.18}), using again (\ref{2.13}) and $\beta
>e^{-\lambda \varepsilon }$, we obtain for sufficiently large $\lambda $%
\begin{equation*}
J_{\lambda ,\beta }\left( u_{1}+h\right) -J_{\lambda ,\beta }\left(
u_{1}\right) -J_{\lambda ,\beta }^{\prime }\left( u_{1}\right) \left(
h\right)
\end{equation*}%
\begin{equation*}
\geq C_{2}e^{-2\lambda \left( c+\varepsilon \right)
}\dint\limits_{G_{c}}\left( \lambda \left( \nabla h\right) ^{2}+\lambda
^{3}h^{2}\right) \varphi _{\lambda }^{2}dx-C_{2}e^{-2\lambda \varepsilon
}\lambda ^{3}\left\Vert h\right\Vert _{H^{k_{n}}\left( G_{c}\right)
}^{2}+\beta \left\Vert h\right\Vert _{H^{k_{n}}\left( G_{c}\right) }^{2}
\end{equation*}%
\begin{equation*}
\geq C_{2}e^{2\lambda \varepsilon }\left\Vert h\right\Vert _{H^{1}\left(
G_{c+2\varepsilon }\right) }^{2}-C_{2}e^{-2\lambda \varepsilon }\lambda
^{3}\left\Vert h\right\Vert _{H^{k_{n}}\left( G_{c}\right) }^{2}+\frac{%
e^{-\lambda \varepsilon }}{2}\left\Vert h\right\Vert _{H^{k_{n}}\left(
G_{c}\right) }^{2}+\frac{\beta }{2}\left\Vert h\right\Vert _{H^{k_{n}}\left(
G_{c}\right) }^{2}
\end{equation*}%
\begin{equation*}
\geq C_{2}e^{2\lambda \varepsilon }\left\Vert h\right\Vert _{H^{1}\left(
G_{c+2\varepsilon }\right) }^{2}+\frac{\beta }{2}\left\Vert h\right\Vert
_{H^{k_{n}}\left( G_{c}\right) }^{2}.\text{ \ \ \ \ }\square
\end{equation*}

\section{Cauchy problem for the quasilinear elliptic equation}

\label{sec:3}

In this section we apply Theorem 2.1 to the Cauchy problem for the
quasilinear elliptic equation. We now rewrite the operator $A$ in (\ref%
{2.100}) as%
\begin{equation}
Au:=L_{ell}\left( u\right) =\dsum\limits_{i,j=1}^{n}a_{i,j}\left( x\right)
u_{x_{i}x_{j}}+A_{1}\left( x,\nabla u,u\right) ,x\in G,  \label{3.1}
\end{equation}%
\begin{equation}
A_{0}u:=L_{0}u=\dsum\limits_{i,j=1}^{n}a_{i,j}\left( x\right) u_{x_{i}x_{j}},
\label{3.2}
\end{equation}%
where $a_{i,j}\left( x\right) =a_{j,i}\left( x\right) ,\forall i,j$ and $%
L_{0}$ is the principal part of the operator $L.$ We impose assumption (\ref%
{2.1002}). The ellipticity of the operator $L_{0}$ means that there exist
two constants $\mu _{1},\mu _{2}>0,\mu _{1}\leq \mu _{2}$ such that%
\begin{equation}
\mu _{1}\left\vert \eta \right\vert ^{2}\leq
\dsum\limits_{i,j=1}^{n}a_{i,j}\left( x\right) \eta _{i}\eta _{j}\leq \mu
_{2}\left\vert \eta \right\vert ^{2},\forall x\in \overline{G},\forall \eta
=\left( \eta _{1},...\eta _{n}\right) \in \mathbb{R}^{n}.  \label{3.4}
\end{equation}%
As above, let $\Gamma \subset \partial G$ be the part of the boundary $%
\partial G$, where the Cauchy data are given. Assume that the equation of $%
\Gamma $ is $\Gamma =\left\{ x\in \mathbb{R}^{n}:x_{1}=p\left( \overline{x}%
\right) ,\overline{x}=\left( x_{2},...,x_{n}\right) \in \Gamma ^{\prime
}\subset \mathbb{R}^{n-1}\right\} $ and that the function $g\in C^{2}\left( 
\overline{\Gamma }^{\prime }\right) .$ Here $\Gamma ^{\prime }\subset 
\mathbb{R}^{n-1}$ is a bounded domain. Changing variables as $x=\left( x_{1},%
\overline{x}\right) \Leftrightarrow \left( x_{1}^{\prime },\overline{x}%
\right) ,$ where $x_{1}^{\prime }=x_{1}-p\left( \overline{x}\right) ,$ we
obtain that in new variables

$\Gamma =\left\{ x\in \mathbb{R}^{n}:x_{1}=0,\overline{x}\in \Gamma ^{\prime
}\right\} .$ Here we keep the same notation for $x_{1}$ as before: for the
simplicity of notations. This change of variables does not affect the
property of the ellipticity of the operator $L$. Let $X>0$ be a certain
number. Thus, without any loss of generality, we assume that 
\begin{equation}
G\subset \left\{ x_{1}>0\right\} ,\text{ }\Gamma =\left\{ x\in \mathbb{R}%
^{n}:x_{1}=0,\left\vert \overline{x}\right\vert <X\right\} \subset \partial
G.  \label{3.5}
\end{equation}

\textbf{Cauchy Problem for the Quasilinear Elliptic Equation}. \emph{Suppose
that conditions (\ref{3.5}) hold. Find such a function }$u\in
H^{k_{n}}\left( G\right) $\emph{\ that satisfies the equation }%
\begin{equation}
L_{ell}\left( u\right) =0  \label{3.6}
\end{equation}%
\emph{and has the following Cauchy data }$g_{0},g_{1}$ \emph{at }$\Gamma $%
\begin{equation}
u\mid _{\Gamma }=g_{0}\left( \overline{x}\right) ,u_{x_{1}}\mid _{\Gamma
}=g_{1}\left( \overline{x}\right) .  \label{3.7}
\end{equation}

Let $\lambda >1$ and $\nu >1$ be two large parameters, which we define
later. Consider two arbitrary numbers $a,c=const.\in \left( 0,1/2\right) ,$
where $a<c$. To introduce the Carleman estimate, consider functions $\psi
\left( x\right) $, $\varphi _{\lambda }\left( x\right) $ defined as%
\begin{equation}
\psi \left( x\right) =x_{1}+\frac{\left\vert \overline{x}\right\vert ^{2}}{%
X^{2}}+a,\varphi _{\lambda }\left( x\right) =\exp \left( \lambda \psi ^{-\nu
}\right) .  \label{3.70}
\end{equation}%
Then the analogs of sets (\ref{2.0})\emph{\ }and\emph{\ }$\Gamma _{c}$ are 
\begin{eqnarray}
G_{c} &=&\left\{ x:x_{1}>0,\left( x_{1}+\frac{\left\vert \overline{x}%
\right\vert ^{2}}{X^{2}}+a\right) ^{-\nu }>c^{-\nu }\right\} ,  \label{3.8}
\\
\xi _{c} &=&\left\{ x:x_{1}>0,\left( x_{1}+\frac{\left\vert \overline{x}%
\right\vert ^{2}}{X^{2}}+a\right) ^{-\nu }=c^{-\nu }\right\} ,  \label{3.80}
\\
\Gamma _{c} &=&\left\{ x:x_{1}=0,\left( \frac{\left\vert \overline{x}%
\right\vert ^{2}}{X^{2}}+a\right) ^{-\nu }>c^{-\nu }\right\} .  \label{3.9}
\end{eqnarray}%
Hence, $\partial G_{c}=\xi _{c}\cup \Gamma _{c}.$ Below in this sections we
keep notations (\ref{3.8})-(\ref{3.9}). We assume that $G_{c}\neq
\varnothing $ and $\overline{G}_{c}\subseteq G.$ By (\ref{3.5}) and (\ref%
{3.9}) $\Gamma _{c}\subset \Gamma .$ For a sufficiently small number $%
\varepsilon >0$ and for $k=1,2$ define the subdomain $G_{c+2\varepsilon }$
as 
\begin{equation}
G_{c+2\varepsilon }=\left\{ x:x_{1}>0,\left( x_{1}+\frac{\left\vert 
\overline{x}\right\vert ^{2}}{X^{2}}+a\right) ^{-\nu }>c^{-\nu
}+2\varepsilon \right\} .  \label{3.91}
\end{equation}%
Hence, $G_{c+2\varepsilon }\subset G_{c}.$ Lemma 3.1 follows immediately
from Lemma 3 of \S 1 of Chapter 4 of the book \cite{LRS}.

\textbf{Lemma 3.1 }(Carleman estimate). \emph{There exist a sufficiently
large number }

$\nu _{0}=\nu _{0}\left( a,\mu _{1},\mu _{2},\max_{i,j}\left\Vert
a_{i,j}\right\Vert _{C^{1}\left( \overline{\Omega }_{c}\right) },X,n\right)
>1$\emph{\ and a sufficiently large absolute constant }$\lambda _{0}>1$\emph{%
\ such that for all }$\nu \geq \nu _{0},\lambda \geq \lambda _{0}$\emph{\
and for all functions }$u\in C^{2}\left( \overline{G}_{1/2}\right) $\emph{\
the following pointwise Carleman estimate is valid for all }$x\in G_{1/2}$%
\emph{\ with a constant }$C=C\left( n,\max_{i,j}\left\Vert
a_{i,j}\right\Vert _{C^{1}\left( \overline{G}_{1/2}\right) }\right) >0$ 
\emph{and with the function }$\varphi _{\lambda }$\emph{\ from (\ref{3.70})} 
\begin{eqnarray*}
\left( L_{0}u\right) ^{2}\varphi _{\lambda }^{2} &\geq &C\lambda \left\vert
\nabla u\right\vert ^{2}\varphi _{\lambda }^{2}+C\lambda ^{3}u^{2}\varphi
_{\lambda }^{2}+\func{div}U, \\
\left\vert U\right\vert &\leq &C\lambda ^{3}\left[ \left( \nabla u\right)
^{2}+u^{2}\right] \varphi _{\lambda }^{2}.
\end{eqnarray*}

This Carleman estimate allows us to construct the weighted Tikhonov
functional to solve the Cauchy problem (\ref{3.6}), (\ref{3.7}).\emph{\ }%
Similarly with (\ref{2.1005}) we minimize the functional $J_{\lambda ,\beta
,ell}\left( u\right) $ (\ref{3.10}) on the set $B\left( R\right) $ defined
in (\ref{1}), where 
\begin{equation}
J_{\lambda ,\beta ,ell}\left( u\right) =e^{-2\lambda \left( c^{-\nu
}+\varepsilon \right) }\dint\limits_{G_{c}}\left[ L_{ell}\left( u\right) %
\right] ^{2}\varphi _{\lambda }^{2}dx+\beta \left\Vert u\right\Vert
_{H^{k_{n}}\left( G_{c}\right) }^{2},  \label{3.10}
\end{equation}%
where the CWF $\varphi _{\lambda }\left( x\right) $ is the one in (\ref{3.70}%
) and functions $g_{0}$ and $g_{1}$ are ones in (\ref{3.7}). Hence, Lemma
3.1 and Theorem 2.1 immediately imply Theorem 3.1.

\textbf{Theorem 3.1.} \emph{Let }$R>0$\emph{\ be an arbitrary number. Let }$%
B\left( R\right) $\emph{\ and }$H_{0}^{k_{n}}\left( G_{c}\right) $\emph{\ be
the sets defined in (\ref{1}) and (\ref{2}) respectively. Then for every
point }$u\in B\left( R\right) $\emph{\ there exists the Fr\'{e}chet
derivative }$J_{\lambda ,\beta ,ell}^{\prime }\left( u\right) \in
H_{0}^{k_{n}}\left( G_{c}\right) .$\emph{\ Choose the numbers }$\nu =\nu
_{0},$ $\lambda _{0}$ \emph{as in Lemma 3.1. There exists a sufficiently
large number }$\lambda _{1}=\lambda _{1}\left( R,L\right) >\lambda _{0}>1$%
\emph{\ such that for all }$\lambda \geq \lambda _{1}$\emph{\ and for every }%
$\beta \in \left( e^{-\lambda \varepsilon },1\right) $\emph{\ the functional 
}$J_{\lambda ,\beta ,ell}\left( u\right) $\emph{\ is strictly convex on the
set }$B\left( R\right) ,$%
\begin{equation*}
J_{\lambda ,\beta ,ell}\left( u_{2}\right) -J_{\lambda ,\beta ,ell}\left(
u_{1}\right) -J_{\lambda ,\beta ,ell}^{\prime }\left( u_{1}\right) \left(
u_{2}-u_{1}\right)
\end{equation*}%
\begin{equation*}
\geq C_{2}e^{2\lambda \varepsilon }\left\Vert u_{2}-u_{1}\right\Vert
_{H^{1}\left( G_{c+2\varepsilon }\right) }^{2}+\frac{\beta }{2}\left\Vert
u_{2}-u_{1}\right\Vert _{H^{k_{n}}\left( G_{c}\right) }^{2},\forall
u_{1},u_{2}\in B\left( R\right) .
\end{equation*}%
\emph{Here the number }$C_{2}=C_{2}\left( L_{ell},R,c\right) >0$\emph{\
depends only on listed parameters. }

\section{Quasilinear parabolic equation with the lateral Cauchy data}

\label{sec:4}

Choose an arbitrary $T=const.>0$ and denote $G_{T}=G\times \left(
-T,T\right) .$ Let $L_{par}$ be the quasilinear elliptic operator of the
second order in $G_{T},$ which we define the same way as the operator $%
L_{ell}$ in (\ref{3.1})-(\ref{3.4}) with the only difference that now its
coefficients depend on both $x$ and $t$ and the domain $G$ is replaced with
the domain $G_{T}.$ Let $L_{0,par}$ be the similarly defined principal part
of the operator $L_{par},$ see (\ref{3.2}). Next, we define the quasilinear
parabolic operator as $P=\partial _{t}-L_{par}$. The principal part of $P$
is $P_{0}=\partial _{t}-L_{0,par}.$ Thus,%
\begin{equation}
L_{par}\left( u\right) =\dsum\limits_{i,j=1}^{n}a_{i,j}\left( x,t\right)
u_{x_{i}x_{j}}+A_{1}\left( x,t,\nabla u,u\right) ,  \label{5.1}
\end{equation}%
\begin{equation}
Au:=Pu=u_{t}-L_{par}\left( u\right) ,\left( x,t\right) \in G_{T},
\label{5.2}
\end{equation}%
\begin{equation}
P_{0}u=u_{t}-L_{0,par}u=u_{t}-\dsum\limits_{i,j=1}^{n}a_{i,j}\left(
x,t\right) u_{x_{i}x_{j}},  \label{5.3}
\end{equation}%
\begin{equation}
a_{i,j}\in C^{1}\left( \overline{G}_{T}\right) ,\text{ }A_{1}\in C\left( 
\overline{G}_{T}\right) \times C^{3}\left( \mathbb{R}^{n+1}\right) ,
\label{5.4}
\end{equation}%
\begin{equation}
\mu _{1}\left\vert \eta \right\vert ^{2}\leq
\dsum\limits_{i,j=1}^{n}a_{i,j}\left( x,t\right) \eta _{i}\eta _{j}\leq \mu
_{2}\left\vert \eta \right\vert ^{2},\forall \left( x,t\right) \in \overline{%
G}_{T},\forall \eta =\left( \eta _{1},...\eta _{n}\right) \in \mathbb{R}^{n}.
\label{5.5}
\end{equation}%
Let $\Gamma \subset \partial G$, $\Gamma \in C^{2}$ be the subsurface of the
boundary $\partial G$ with the same properties as ones in section 3.\ Denote 
$\Gamma _{T}=\Gamma \times \left( -T,T\right) .$ Without loss of generality
we assume that $\Gamma $ is the same as in (\ref{3.5}). Consider the
parabolic equation 
\begin{equation}
P\left( u\right) =u_{t}-L_{par}\left( u\right) =0\text{ \ in }G_{T}.\text{ }
\label{4.2}
\end{equation}

\textbf{Cauchy Problem with the Lateral Data for Quasilinear Parabolic
Equation (\ref{4.2}).} \emph{Assume that conditions (\ref{3.5}) hold. Find
such a function }$u\in H^{k_{n+1}}\left( G_{T}\right) $\emph{\ that
satisfies equation (\ref{4.2}) and has the following lateral Cauchy data }$%
g_{0},g_{1}$ \emph{at }$\Gamma _{T}$%
\begin{equation}
u\mid _{\Gamma _{T}}=g_{0}\left( \overline{x},t\right) ,u_{x_{1}}\mid
_{\Gamma _{T}}=g_{1}\left( \overline{x},t\right) .  \label{4.3}
\end{equation}

We now introduce the Carleman estimate which is similar with the one of
section 3. Let $\lambda >1$ and $\nu >1$ be two large parameters, which we
define later. Consider two arbitrary numbers $a,c=const.\in \left(
0,1/2\right) ,$ where $a<c$. Consider functions $\psi \left( x,t\right) $, $%
\varphi _{\lambda }\left( x,t\right) $ defined as%
\begin{equation}
\psi \left( x,t\right) =x_{1}+\frac{\left\vert \overline{x}\right\vert ^{2}}{%
X^{2}}+\frac{t^{2}}{T^{2}}+a,\text{ }\varphi _{\lambda }\left( x,t\right)
=\exp \left( \lambda \psi ^{-\nu }\right) .  \label{4.30}
\end{equation}%
Analogs of conditions (\ref{3.8})-(\ref{3.91})\emph{\ }are 
\begin{eqnarray}
G_{T,c} &=&\left\{ \left( x,t\right) :x_{1}>0,\left( x_{1}+\frac{\left\vert 
\overline{x}\right\vert ^{2}}{X^{2}}+\frac{t^{2}}{T^{2}}+a\right) ^{-\nu
}>c^{-\nu }\right\} ,  \label{4.4} \\
\xi _{c} &=&\left\{ \left( x,t\right) :x_{1}>0,\left( x_{1}+\frac{\left\vert 
\overline{x}\right\vert ^{2}}{X^{2}}+\frac{t^{2}}{T^{2}}+a\right) ^{-\nu
}=c^{-\nu }\right\} ,  \label{4.40} \\
\Gamma _{c} &=&\left\{ \left( x,t\right) :x_{1}=0,\left( \frac{\left\vert 
\overline{x}\right\vert ^{2}}{X^{2}}+\frac{t^{2}}{T^{2}}+a\right) ^{-\nu
}>c^{-\nu }\right\} ,  \label{4.5} \\
\partial G_{T,c} &=&\xi _{c}\cup \Gamma _{c},  \label{4.6} \\
G_{c+2\varepsilon ,T} &=&\left\{ \left( x,t\right) :x_{1}>0,\left( x_{1}+%
\frac{\left\vert \overline{x}\right\vert ^{2}}{X^{2}}+\frac{t^{2}}{T^{2}}%
+a\right) ^{-\nu }>c^{-\nu }+2\varepsilon \right\} .  \label{4.7}
\end{eqnarray}%
We assume that%
\begin{equation}
G_{T,c}\neq \varnothing ,G_{T,c}\subset G_{T}\text{ and }\overline{G}%
_{T,c}\cap \left\{ t=\pm T\right\} =\varnothing .  \label{4.700}
\end{equation}%
In (\ref{4.7}) $\varepsilon >0$ is so small that $G_{c+2\varepsilon ,T}\neq
\varnothing .$ Below in this section we use notations (\ref{4.30})-(\ref{4.7}%
). By (\ref{4.1}) and (\ref{4.5}) $\Gamma _{c}\subset \Theta _{T}.$ Lemma
4.1 follows immediately from Lemma 3 of \S 1 of chapter 4 of the book \cite%
{LRS}.

\textbf{Lemma 4.1 }(Carleman estimate). \emph{Let }$P_{0}$\emph{\ be the
parabolic operator defined via (\ref{5.1})-(\ref{5.5}). There exist a
sufficiently large number }$\nu _{0}=\nu _{0}\left( a,c,\mu _{1},\mu
_{2},\max_{i,j}\left\Vert a_{i,j}\right\Vert _{C^{1}\left( \overline{\Omega }%
_{c}\right) },X,T\right) >1$\emph{\ and a sufficiently large absolute
constant }$\lambda _{0}>1$\emph{\ such that for all }$\nu \geq \nu
_{0},\lambda \geq \lambda _{0}$\emph{\ and for all functions }$u\in
C^{2,1}\left( \overline{G}_{T,1/2}\right) $\emph{\ the following pointwise
Carleman estimate is valid for all }$\left( x,t\right) \in G_{T,1/2}$\emph{\
with a constant }$C=C\left( n,\max_{i,j}\left\Vert a_{i,j}\right\Vert
_{C^{1}\left( \overline{G}\right) }\right) >0$ \emph{and} \emph{with the
function }$\varphi _{\lambda }$\emph{\ defined in (\ref{4.30})} 
\begin{eqnarray*}
\left( P_{0}u\right) ^{2}\varphi _{\lambda }^{2} &\geq &C\lambda \left\vert
\nabla u\right\vert ^{2}\varphi _{\lambda }^{2}+C\lambda ^{3}u^{2}\varphi
_{\lambda }^{2}+\func{div}U+V_{t}, \\
\left\vert U\right\vert ,\left\vert V\right\vert &\leq &C\lambda ^{3}\left[
\left( \nabla u\right) ^{2}+u_{t}^{2}+u^{2}\right] \varphi _{\lambda }^{2}.
\end{eqnarray*}

Let $R>0$ be an arbitrary number. Similarly with (\ref{1}) and (\ref{2}) let 
\begin{equation}
B\left( R\right) =\left\{ u\in H^{k_{n+1}}\left( G_{T,c}\right) :\left\Vert
u\right\Vert _{H^{k_{n+1}}\left( G_{T,c}\right) }<R,u\mid _{\Gamma
_{c}}=g_{0}\left( \overline{x},t\right) ,\partial _{n}u\mid _{\Gamma
_{c}}=g_{1}\left( \overline{x},t\right) \right\} ,  \label{5.6}
\end{equation}%
\begin{equation}
H_{0}^{k_{n+1}}\left( G_{T,c}\right) =\left\{ u\in H^{k_{n+1}}\left(
G_{T,c}\right) :u\mid _{\Gamma _{c}}=\partial _{n}u\mid _{\Gamma
_{c}}=0\right\} .  \label{5.7}
\end{equation}

To solve the problem (\ref{4.2}), (\ref{4.3}), we minimize the functional $%
J_{\lambda ,\beta ,par}\left( u\right) $ in (\ref{4.8}) on the set $B\left(
R\right) $ defined (\ref{5.6}), where 
\begin{equation}
J_{\lambda ,\beta ,par}\left( u\right) =e^{-2\lambda \left( c^{-\nu
}+\varepsilon \right) }\dint\limits_{G_{T,c}}\left[ P\left( u\right) \right]
^{2}\varphi _{\lambda }^{2}dx+\beta \left\Vert u\right\Vert
_{H^{k_{n+1}}\left( G_{T,c}\right) }^{2},  \label{4.8}
\end{equation}%
where the operator $P$ is defined via (\ref{5.1})-(\ref{5.5}). Hence, Lemma
4.1 and Theorem 2.1 imply Theorem 4.1.

\textbf{Theorem 4.1. }\emph{Let }$R>0$\emph{\ be an arbitrary number.\ Let }$%
B\left( R\right) $\emph{\ and }$H_{0}^{k_{n+1}}\left( G_{T,c}\right) $ \emph{%
be the sets defined in (\ref{5.6}) and (\ref{5.7}) respectively. Then} \emph{%
for every point }$u\in B\left( R\right) $\emph{\ there exists the Fr\'{e}%
chet derivative }$J_{\lambda ,\beta ,par}^{\prime }\left( u\right) \in
H_{0}^{k_{n+1}}\left( G_{T,c}\right) .$\emph{\ Choose the numbers }$\nu =\nu
_{0},$ $\lambda _{0}$ \emph{as in Lemma 4.1. There exists a sufficiently
large number }$\lambda _{1}=\lambda _{1}\left( R,P\right) >\lambda _{0}>1$%
\emph{\ such that for all }$\lambda \geq \lambda _{1}$\emph{\ and for every }%
$\beta \in \left( e^{-\lambda \varepsilon },1\right) $\emph{\ the functional 
}$J_{\lambda ,\beta ,par}\left( u\right) $\emph{\ is strictly convex on the
set }$B\left( R\right) ,$%
\begin{equation*}
J_{\lambda ,\beta ,par}\left( u_{2}\right) -J_{\lambda ,\beta ,par}\left(
u_{1}\right) -J_{\lambda ,\beta ,par}^{\prime }\left( u_{1}\right) \left(
u_{2}-u_{1}\right)
\end{equation*}%
\begin{equation*}
\geq C_{3}e^{2\lambda \varepsilon }\left\Vert u_{2}-u_{1}\right\Vert
_{H^{1}\left( G_{T,c+2\varepsilon }\right) }^{2}+\frac{\beta }{2}\left\Vert
u_{2}-u_{1}\right\Vert _{H^{k_{n}+1}\left( G_{T,c}\right) }^{2},\forall
u_{1},u_{2}\in B\left( R\right) .
\end{equation*}%
\emph{Here the number }$C_{3}=C_{3}\left( L_{par},R,c\right) >0$\emph{\
depends only on listed parameters. }

\section{Quasilinear hyperbolic equation with lateral Cauchy data}

\label{sec:5}

In this section, notations for the domain $G\subset \mathbb{R}^{n}$ and the
time cylinder $G_{T}$ are the same as ones in section 4. Denote $%
S_{T}=\partial G\times \left( -T,T\right) .$ The Carleman estimate of Lemma
5.1 of this section for was proved in Theorem 1.10.2 of the book of Beilina
and Klibanov \cite{BK}. Other forms of Carleman estimates for the hyperbolic
case can be found in, e.g. Theorem 3.4.1 of the book of Isakov \cite{Is},
Theorem 2.2.4 of the book of Klibanov and Timonov \cite{KT}, Lemma 2 of \S 4
of chapter 4 of the book of Lavrentiev, Romanov and Shishatskii \cite{LRS}
and in Lemma 3.1 of Triggiani and Yao \cite{Trig}.

Let numbers $a_{l},a_{u}>0$ and $a_{l}<a_{u}.$ For $x\in G,$ let the
function $a\left( x\right) $ satisfy the following conditions in $G$%
\begin{equation}
a\left( x\right) \in \left[ a_{l},a_{u}\right] ,a\in C^{1}\left( \overline{G}%
\right) .  \label{6.1}
\end{equation}%
In addition, we assume that there exists a point $x_{0}\in G$ such that%
\begin{equation}
\left( \nabla a\left( x\right) ,x-x_{0}\right) \geq 0,\forall x\in \overline{%
G},  \label{6.2}
\end{equation}%
where $\left( \cdot ,\cdot \right) $ denotes the scalar product in $\mathbb{R%
}^{n}$.\emph{\ }In particular, if $a\left( x\right) \equiv 1,$ then (\ref%
{6.2}) holds for any $x_{0}\in G.$ We need inequality (\ref{6.2}) for the
validity of the Carleman estimate of Lemma 5.1. Assume that the function $%
A_{1}$ satisfies condition (\ref{5.4}). Consider the quasilinear hyperbolic
equation in the time cylinder $G_{T}$ with the lateral Cauchy data $%
g_{0}\left( x,t\right) ,g_{1}\left( x,t\right) ,$ 
\begin{eqnarray}
L_{hyp}u &=&a\left( x\right) u_{tt}-\Delta u-A_{1}\left( x,t,\nabla
u,u\right) =0\text{ in }G_{T},  \label{6.3} \\
u &\mid &_{S_{T}}=g_{0}\left( x,t\right) ,\partial _{n}u\mid
_{S_{T}}=g_{1}\left( x,t\right) .  \label{6.4}
\end{eqnarray}%
Denote $L_{0,hyp}u=a\left( x\right) u_{tt}-\Delta u.$

\textbf{Cauchy Problem with the Lateral Data for the Hyperbolic Equation (%
\ref{6.3})}. \emph{Find the function }$u\in H^{k_{n+1}}\left( G_{T}\right) $%
\emph{\ satisfying conditions (\ref{6.3}), (\ref{6.4}).}

Let the number $\eta \in \left( 0,1\right) .$ Let $\lambda >1$ be a large
parameter$.$ Define functions $\xi \left( x,t\right) $ and $\varphi
_{\lambda }\left( x,t\right) $ as%
\begin{equation}
\xi \left( x,t\right) =\left\vert x-x_{0}\right\vert ^{2}-\eta t^{2},\varphi
_{\lambda }\left( x,t\right) =\exp \left[ \lambda \psi \left( x,t\right) %
\right] .  \label{6.5}
\end{equation}%
Similarly with (\ref{2.0}), for a number $c>0$ define the hypersurface $\xi
_{c}$ and the domain $G_{T,c}$ as%
\begin{equation}
\xi _{c}=\left\{ \left( x,t\right) \in G_{T}:\xi \left( x,t\right)
=c\right\} ,\text{ }G_{T,c}=\left\{ \left( x,t\right) \in G_{T}:\xi \left(
x,t\right) >c\right\} .  \label{6.6}
\end{equation}

\textbf{Lemma 5.1 }(Carleman estimate). \emph{Let }$n\geq 2$\emph{\ and
conditions (\ref{6.1}) be satisfied. Also, assume that there exists a point }%
$x_{0}\in G$\emph{\ such that (\ref{6.2}) holds. Denote }$M=M\left(
x_{0},G\right) =\max_{x\in \overline{G}}\left\vert x-x_{0}\right\vert .$%
\emph{\ Choose such a number }$c>0$\emph{\ that }$G_{T,c}\neq \varnothing .$%
\emph{\ Let }$\varphi _{\lambda }\left( x,t\right) $\emph{\ be the function
defined in (\ref{6.5}), sets }$\xi _{c},G_{T,c}$\emph{\ be the ones defined
in (\ref{6.6}) and conditions (\ref{4.700}) are valid for the case of the
domain }$G$\emph{\ of this section. Then there exists a number }$\eta
_{0}=\eta _{0}\left( G,M,a_{l},a_{u},\left\Vert \nabla a\right\Vert
_{C\left( \overline{G}\right) }\right) \in \left( 0,1\right) $\emph{\ such
that for any }$\eta \in \left( 0,\eta _{0}\right) $\emph{\ one can choose a
sufficiently large number }$\lambda _{0}=\lambda _{0}\left(
G,M,a_{l},a_{u},\left\Vert \nabla a\right\Vert _{C\left( \overline{G}\right)
},\eta _{0},c\right) >1$\emph{\ and the number }$C_{4}=C_{4}\left(
G,M,a_{l},a_{u},\left\Vert \nabla a\right\Vert _{C\left( \overline{G}\right)
},\eta _{0},c\right) >0$\emph{, such that for all }$u\in C^{2}\left( 
\overline{G}_{T,c}\right) $\emph{\ and for all }$\lambda \geq \lambda _{0}$%
\emph{\ the following pointwise Carleman estimate is valid} 
\begin{equation*}
\left( L_{0,hyp}u\right) ^{2}\varphi _{\lambda }^{2}\geq C_{4}\lambda \left(
\left\vert \nabla u\right\vert ^{2}+u_{t}^{2}\right) \varphi _{\lambda
}^{2}+C_{4}\lambda ^{3}u^{2}\varphi _{\lambda }^{2}+\func{div}U+V_{t}\text{
in }G_{T,c},
\end{equation*}%
\emph{\ } 
\begin{equation*}
\left\vert U\right\vert ,\left\vert V\right\vert \leq C_{4}\lambda
^{3}\left( \left\vert \nabla u\right\vert ^{2}+u_{t}^{2}+u^{2}\right)
\varphi _{\lambda }^{2}.
\end{equation*}%
\emph{\ } \emph{In the case }$a\left( x\right) \equiv 1$\emph{\ one }$\eta
_{0}=1.$

Again, let $R>0$ be an arbitrary number. Similarly with (\ref{5.6}) and (\ref%
{5.7}) let%
\begin{equation}
B\left( R\right) =\left\{ u\in H^{k_{n+1}}\left( G_{T,c}\right) :\left\Vert
u\right\Vert _{H^{k_{n+1}}\left( G_{T,c}\right) }<R,u\mid
_{S_{T}}=g_{0}\left( x,t\right) ,\partial _{n}u\mid _{S_{T}}=g_{1}\left(
x,t\right) \right\} ,  \label{6.9}
\end{equation}
\begin{equation}
H_{0}^{k_{n+1}}\left( G_{T,c}\right) =\left\{ u\in H^{k_{n+1}}\left(
G_{T,c}\right) :u\mid _{S_{T}}=\partial _{n}u\mid _{S_{T}}=0\right\} .
\label{6.10}
\end{equation}%
To solve the Cauchy problem posed in this section, we minimize the
functional $J_{\lambda ,\beta ,hyp}\left( u\right) $ in (\ref{6.7}) on the
set $B\left( R\right) $ defined in (\ref{6.9}), where 
\begin{equation}
J_{\lambda ,\beta ,hyp}\left( u\right) =e^{-2\lambda \left( c+\varepsilon
\right) }\dint\limits_{G_{T,c}}\left[ L_{hyp}\left( u\right) \right]
^{2}\varphi _{\lambda }^{2}dx+\beta \left\Vert u\right\Vert
_{H^{k_{n+1}}\left( G_{T,c}\right) }^{2}.  \label{6.7}
\end{equation}%
Hence, Lemma 5.1 and Theorem 2.1 imply Theorem 5.1.

\textbf{Theorem 5.1}. \emph{Let }$R>0$\emph{\ be an arbitrary number.\ Let }$%
B\left( R\right) $\emph{\ and }$H_{0}^{k_{n+1}}\left( G_{T,c}\right) $ \emph{%
be the sets defined in (\ref{6.9}) and (\ref{6.10}) respectively. Then} 
\emph{for every point }$u\in B\left( R\right) $\emph{\ there exists the Fr%
\'{e}chet derivative }$J_{\lambda ,\beta ,hyp}^{\prime }\left( u\right) \in
H_{0}^{k_{n+1}}\left( G_{T,c}\right) .$\emph{\ Let} $\lambda _{0}=\lambda
_{0}\left( G,M,a_{l},a_{u},\left\Vert \nabla a\right\Vert _{C\left( 
\overline{G}\right) },\eta _{0},c\right) >1$ \emph{be the number of Lemma
5.1. There exists a sufficiently large number }

$\lambda _{1}=\lambda _{1}\left( R,L_{hyp},G,M,a_{l},a_{u},\left\Vert \nabla
a\right\Vert _{C\left( \overline{G}\right) },\eta _{0},c\right) >\lambda
_{0}>1$\emph{\ such that for all }$\lambda \geq \lambda _{1}$\emph{\ and for
every }$\beta \in \left( e^{-\lambda \varepsilon },1\right) $\emph{\ the
functional }$J_{\lambda ,\beta ,hyp}\left( u\right) $\emph{\ is strictly
convex on the set }$B\left( R\right) ,$%
\begin{equation*}
J_{\lambda ,\beta ,par}\left( u_{2}\right) -J_{\lambda ,\beta ,hyp}\left(
u_{1}\right) -J_{\lambda ,\beta ,hyp}^{\prime }\left( u_{1}\right) \left(
u_{2}-u_{1}\right)
\end{equation*}%
\begin{equation*}
\geq C_{5}e^{2\lambda \varepsilon }\left\Vert u_{2}-u_{1}\right\Vert
_{H^{1}\left( G_{T,c+2\varepsilon }\right) }^{2}+\frac{\beta }{2}\left\Vert
u_{2}-u_{1}\right\Vert _{H^{k_{n}+1}\left( G_{T,c}\right) }^{2},\forall
u_{1},u_{2}\in B\left( R\right) .
\end{equation*}%
\emph{Here the number }$C_{5}=C_{3}\left( L_{hyp},R,c\right) >0$\emph{.}

\end{document}